\documentclass[11pt,a4paper,1p,sort&compress]{elsarticle}
\usepackage[]{times}
\usepackage{fullpage}
\usepackage{graphicx}
\usepackage{subfigure}
\usepackage{amsmath}
\usepackage{fancyhdr}
\date{}%leave empty

\usepackage{bm}
\usepackage{amsfonts,amsmath,amssymb,amsthm}
\usepackage{graphicx}
\usepackage{subfigure}
\usepackage{wrapfig}
\usepackage{float}

\begin{document}

\begin{frontmatter}

\title{Simulation of two-fluid flows using a Finite Element/level set method. Application to
bubbles and vesicle dynamics}

\author[ad1]{V. Doyeux}
\ead{vincent.doyeux@ujf-grenoble.fr}
\author[ad1]{Y. Guyot}
\author[ad2]{V. Chabannes}
\ead{vincent.chabannes@imag.fr}
\author[ad2,ad3]{C. Prud'homme}
\ead{christophe.prudhomme@ujf-grenoble.fr}
\author[ad1]{M. Ismail\corref{cor}}
\ead{mourad.ismail@ujf-grenoble.fr}

\address[ad1]{Universit\'e Grenoble 1 / CNRS, Laboratoire Interdisciplinaire de
  Physique / UMR 5588 Grenoble, F-38041, France}
\address[ad2]{Universit\'e Grenoble 1 / CNRS, Laboratoire Jean Kuntzman / UMR 5224. Grenoble, F-38041, France}
\address[ad3]{Universit\'e de Strasbourg / CNRS, IRMA / UMR  7501. Strasbourg, F-67000, France}

\cortext[cor]{Corresponding author}

\begin{abstract}
  A new framework for two-fluids flow using a Finite Element/Level Set method is presented and
  verified through the simulation of the rising of a bubble in a viscous
  fluid. This model is then enriched to deal with vesicles (which mimic red
  blood cells mechanical behavior) by introducing a Lagrange multiplier to constrain the
  inextensibility of the membrane. Moreover, high order polynomial approximation
  is used to increase the accuracy of the simulations. A validation of this
  model is finally presented on known behaviors of vesicles under flow such as
  ``tank treading'' and tumbling motions.
\end{abstract}
\end{frontmatter}

%\listoftodos

\section*{Introduction}
\label{sec:introdution}

Vesicles are systems of two-fluids separated by a bi-layer membrane of
phospholipids which has the property to be inextensible. These objects are
biomimetics in the sense that they reproduce some biological objects
behaviors. Specifically, vesicles have a mechanical behavior close to the one of
Red Blood Cells (\textbf{RBC}) in a fluid flow. Indeed, it has been accepted for
many years as a good model for RBC and they have been studied experimentally,
theoretically and numerically. Simulating vesicles is very challenging in the
sense that it combines fluid structure interaction and two-fluid flow
systems. Several methods have already been developed such as lattice Boltzmann
methods \cite{Kaoui2011}, boundary integral methods \cite{Beaucourt2004}, or
level set methods using finite difference method \cite{Maitre2010, Salac2011} or
finite element method \cite{Laadhari2011}.
Recently, another model based on a ``Necklace'' of rigid particles was proposed
by one of the authors to model vesicles dynamic in fluid flow \cite{ismail_lefebvre:2012}.

We present in this paper a new framework
to simulate vesicles by level set method using finite element
approximations. This framework has been inspired by \cite{Winkelmann2007} and
\cite{Laadhari2011} albeit with some differences in the strategy (mesh
adaptation, Lagrange multipliers on advection equation). We propose to verify
the framework, from the numerical point of view in a first time --- using a
benchmark for two-fluid flow by level set method which consists of the rising of
a bubble in a viscous fluid. --- Then we present our strategy for the simulation of
vesicles and validate it on some known behaviors of vesicles under flow as the
\textit{tank treading} and \textit{tumbling} motions.

\section{Level set description}
\label{sec:level_set_description}

\subsection{Description}

Let's define a bounded domain $\Omega \subset \mathbb{R}^p$ ($p=2,3$) decomposed
into two subdomains $\Omega_1$ and $\Omega_2$. We denote $\Gamma$ the interface
between the two partitions. The goal of the level set method is to track
implicitly the interface $\Gamma(t)$ moving at a velocity $\bm{u}$. The level
set method has been described in \cite{Osher1988, book_Sethian, book_Osher} and
its main ingredient is a continuous scalar function $\phi$ (the \textit{level
set} function) defined on the whole domain. This function is chosen to be
positive in $\Omega_1$, negative in $\Omega_2$ and zero on $\Gamma$. The motion
of the interface is then described by the advection of the level set function
with a divergence free velocity field $\bm{u}$:
\begin{equation}
  \label{eq:advection}
   \frac{\partial \phi}{\partial t} + \bm{u} \cdot \nabla \phi = 0,\quad \nabla
   \cdot \bm{u} = 0.
\end{equation}
A convenient choice for $\phi$ is a signed distance function to the
interface. Indeed, the property $|\nabla \phi| = 1$ of distance functions  eases the
numerical solution and gives a convenient support for delta and Heaviside
functions (see section \ref{subsec:inter_rel_quant}).
Nevertheless, it is known that the advection equation \eqref{eq:advection} does not conserve
the property $|\nabla \phi|=1$. Thus, when $|\nabla \phi|$ is far from $1$ we use a
fast marching method (\textbf{FMM}) which resets $\phi$ as a distance function
without moving the interface (see \cite{Winkelmann2007} for details about the
fast marching method).

\subsection{Interface related quantities}
\label{subsec:inter_rel_quant}

In two-fluid flow simulations, we need to define some quantities related to the
interface such as the density, the viscosity, or some interface forces. To this end,
we introduce the smoothed Heaviside and delta functions :\\
\small
\begin{minipage}{0.49\textwidth}
  \begin{equation*}
    H_{\varepsilon}(\phi) = \left\{ \begin{array}{cc}
      0, & \phi \leq - \varepsilon,\\
      \displaystyle\frac{1}{2} \left(1+\frac{\phi}{\varepsilon}+\frac{\sin(\frac{\pi \phi}{\varepsilon})}{\pi}\right),  & -\varepsilon \leq \phi \leq \varepsilon, \\
      1, & \phi \geq \varepsilon. \end{array} \right.
  \end{equation*}
\end{minipage}
\begin{minipage}{0.49\textwidth}
  \begin{equation*}
    \delta_{\varepsilon}(\phi) = \left\{ \begin{array}{cc}
      0, & \phi \leq - \varepsilon,\\
      \displaystyle\frac{1}{2 \varepsilon} \left(1+\cos(\frac{\pi \phi}{\varepsilon})\right),  & -\varepsilon \leq \phi \leq \varepsilon, \\
      0, & \phi \geq \varepsilon. \end{array} \right.
  \end{equation*}
\end{minipage}\\

\normalsize
where $\varepsilon$ is a parameter defining a ``numerical thickness'' of the interface. A
typical value of $\varepsilon$ is $1.5 h$ where $h$ is the mesh size of elements
crossed by the iso-value $0$ of the level set function.

The Heaviside function is used to define parameters having different values on
each subdomains. For example, we define the density of two-fluid flow as
$\rho = \rho_2 + (\rho_1-\rho_2) H_{\varepsilon}(\phi)$ (we use a similar
expression for the viscosity $\nu$). Regarding the delta function, it is used to define
quantities on the interface. In particular, in the variational formulations, we
replace integrals over the interface $\Gamma$ by integrals over the entire
domain $\Omega$ using
the smoothed delta function. If $\phi$ is a signed distance function, we have :
$\int_{\Gamma} 1 \simeq \int_{\Omega} \delta_{\varepsilon}(\phi)$.  If $\phi$ is
not close enough to a distance function, then $\int_{\Gamma} 1 \simeq
\int_{\Omega} |\nabla \phi| \delta_{\varepsilon}(\phi)$ which still tends to the
measure of $\Gamma$ as $\varepsilon$ vanishes. However, if $\phi$ is not a
distance function, the support of $\delta_{\varepsilon}$ can have a different
size on each side of the interface. More precisely, the support of
$\delta_{\varepsilon}$ is narrowed on the side where $|\nabla \phi|>1$ and
enlarged on regions where $|\nabla \phi|<1$.  It has been shown in
\cite{Cottet2010} that replacing $\phi$ by $\frac{\phi}{|\nabla \phi|}$ has the
property that $|\nabla \frac{\phi}{|\nabla \phi|}| \simeq 1$ near the interface
and has the same iso-value $0$ as $\phi$. Thus, replacing $\phi$ by
$\frac{\phi}{|\nabla \phi|}$ as support of the delta function does not move the
interface. Moreover, the spread interface has the same size on each part of the
level-set $\phi=0$. It reads $\int_{\Gamma} 1 \simeq \int_{\Omega}
\delta_{\varepsilon}(\frac{\phi}{| \nabla \phi|})$.  The same technique is used for
the Heaviside function.

\subsection{Numerical implementation and coupling with the fluid solver}

We use the finite element C++ library \texttt{Feel++}
\cite{PRUDHOMME:2012:HAL-00662868:3,prudhomme:_feel,Prudhomme.ea:2006} to
discretize and solve the problem.  Equation \eqref{eq:advection} is solved using
a stabilized finite element method. We have implemented several stabilization
methods such as Streamline Upwind Diffusion (SUPG), Galerkin Least Square (GLS)
and Subgrid Scale (SGS). A general review of these methods is available in
\cite{Franca1992}. Other available methods include the Continuous Interior
Penalty method (CIP) are described in \cite{Burman2006}.  The variational
formulation at the semi-discrete level for the stabilized equation \eqref{eq:advection} reads, find $\phi_h
\in {\mathbb R}_h^k$ such that $\forall \psi_h \in {\mathbb R}_h^k$ :
\begin{equation}
  \left(\int_{\Omega} \frac{\partial \phi_h}{\partial t}  \psi_h + \int_{\Omega} (\bm{u}_h \cdot \nabla \phi_h) \psi_h\right) +  S(\phi_h, \psi_h) = 0,
  \label{eq:varadvection}
\end{equation}
where $S(\cdot, \cdot)$ stands for the stabilization bilinear form (see section
\ref{sec:membr-inext} for description of ${\mathbb R}_h^k$ and $\bm{u}_h$).
In our case, we use a Crank-Nicholson scheme which needs only the
solution at previous time step to compute the one at present time. % Indeed, at the

\section{Validation of two-fluid flow solver}

The previous section described the strategy we used to track the interface. We
couple it now to the Navier Stokes equation solver described in
\cite{chabannes11:_high}. In the current section, we present a validation of
this two-fluid flow framework. To do this, we chose to compare our results to
the ones given by the benchmark introduced in
\cite{Hysing2009}.

\subsection{Benchmark problem}

The benchmark objective is to simulate the rise of a 2D bubble in a Newtonian
fluid. The equations solved are the incompressible Navier Stokes equations for
the fluid and the advection for the level set:
\begin{eqnarray}
\rho( \phi(\bm{x}) ) \left(\frac{\partial \bm{u}}{\partial t} + \bm{u} \cdot \nabla \bm{u} \right) + \nabla p - \nabla \cdot \left( \nu(\phi(\bm{x})) (\nabla \bm{u} + (\nabla \bm{u})^T) \right) &=& \rho ( \phi(\bm{x}) ) \bm{g}, \label{eq:bubble_problem1}\\
\nabla \cdot \bm{u} &=& 0, \label{eq:bubble_problem2}\\
\frac{\partial \phi}{\partial t} + \bm{u} \cdot \nabla \phi &=& 0, \label{eq:bubble_problem3}
\end{eqnarray}
where $\rho$ is the density of the fluid, $\nu$ its viscosity, and $\bm{g} \approx (0,
0.98)^T$ is the gravity acceleration.

The computational domain is $\Omega \times ]0, T]$ where $\Omega = (0,1) \times
(0,2)$ and $T=3$. We denote $\Omega_1$ the domain outside the bubble $ \Omega_1
= \{\bm{x} | \phi(\bm{x})>0 \} $, $\Omega_2$ the domain inside the bubble $
\Omega_2 = \{\bm{x} | \phi(\bm{x})<0 \} $ and $\Gamma$ the interface $ \Gamma =
\{\bm{x} | \phi(\bm{x})=0 \} $.  On the lateral walls, slip boundary conditions
are imposed, \emph{i.e.} $\bm{u} \cdot \bm{n} = 0$ and $\bm{t} \cdot (\nabla
\bm{u} + (\nabla \bm{u})^T) \cdot \bm{n}=0$ where $\bm{n}$ is the unit normal to
the interface and $\bm{t}$ the unit tangent. No slip boundary conditions are
imposed on the horizontal walls \emph{i.e.} $\bm{u}=\bm{0}$.  The initial bubble
is circular with a radius $r_0 = 0.25$ and centered on the point $(0.5, 0.5)$.
A surface tension force $\bm{f}_{st}$ is applied on $\Gamma$, it reads :
$\bm{f}_{st} = \int_{\Gamma} \sigma \kappa \bm{n} \simeq \int_{\Omega} \sigma
\kappa \bm{n} \delta_{\varepsilon}(\phi)$ where $\sigma$ stands for the surface
tension between the two-fluids and $\kappa = \nabla \cdot (\frac{\nabla
  \bm{\phi}}{|\nabla \phi|})$ is the curvature of the interface. Note that the
 normal vector $\bm{n}$ is defined here as $\bm{n}=\frac{\nabla
  \phi}{|\nabla \phi|}$.

  We denote with indices $1$ and $2$ the quantities relative to the fluid in
  respectively in $\Omega_1$ and $\Omega_2$. The parameters of the benchmark
  are $\rho_1$, $\rho_2$, $\nu_1$, $\nu_2$ and $\sigma$ and we define two
  dimensionless numbers: first, the Reynolds number which is the
  ratio between inertial and viscous terms and is defined as : $Re = \dfrac{\rho_1 \sqrt{|\bm{g}|
  (2r_0)^3}}{\nu_1}$; second, the E\"otv\"os number which represents the ratio
  between the gravity force and the surface tension $E_0 = \dfrac{4 \rho_1
  |\bm{g}| r_0^2}{\sigma}$.  The table \ref{tab:test_param} reports the values
  of the parameters used for two different test cases proposed in~\cite{Hysing2009}.
\begin{table}[htb!]
  \normalsize
  \begin{center}
    \begin{tabular}{|c|c|c|c|c|c|c|c|} \hline
      Tests & $\rho_1$ & $\rho_2$ & $\nu_1$ & $\nu_2$ & $\sigma$ & $Re$ & $E_0$ \\ \hline
      Test 1 (ellipsoidal bubble) & 1000    & 100      &  10    &   1      &    24.5  & 35   &  10  \\  \hline
      Test 2 (skirted bubble)    &  1000    &  1       &  10    &  0.1     &  1.96    & 35   & 125   \\ \hline
    \end{tabular}
  \end{center}
  \caption{Numerical parameters taken for the benchmarks.}
  \label{tab:test_param}
\end{table}
The quantities measured in \cite{Hysing2009} are $\bm{X_c}$ the center of mass
of the bubble, $\bm{U_c}$ its velocity and the circularity defined as the ratio
between the perimeter of a circle which has the same area and the perimeter
of the bubble which reads $c = \dfrac{2(\pi \int_{\Omega_2} 1 )^{\frac{1}{2}}}{\int_{\Gamma} 1}$.

\subsection{Results}

\begin{figure}[h!tbp]%
  \centering
  \subfigure[Shape at final time ($t=3$).]{
    \includegraphics[scale=0.18]{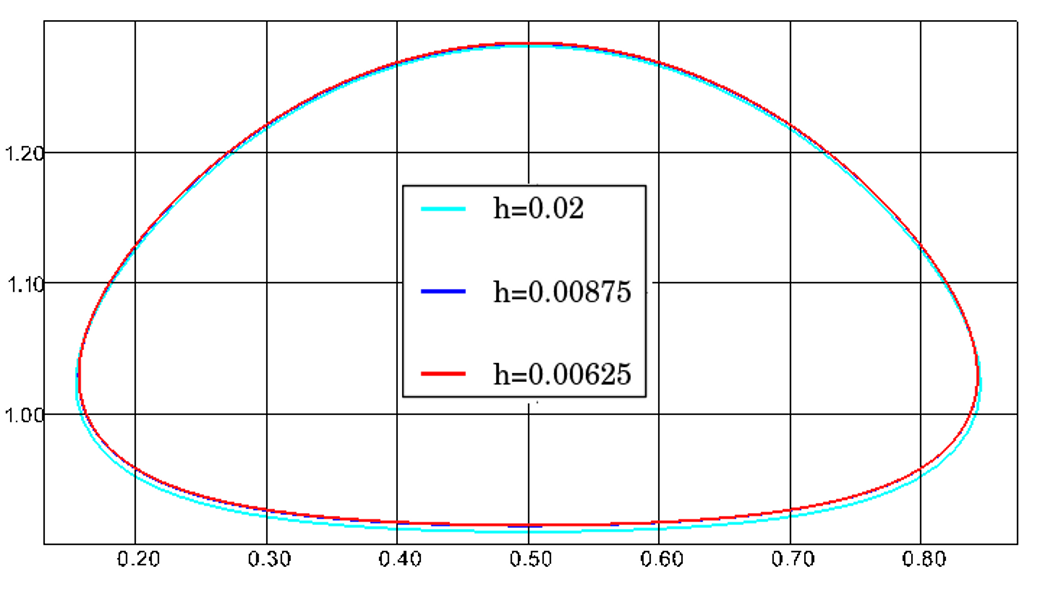}
    \label{subfig:elli_sh}
  }
  \subfigure[$y_c$ vertical position.]{
    \includegraphics[width=0.4\textwidth]{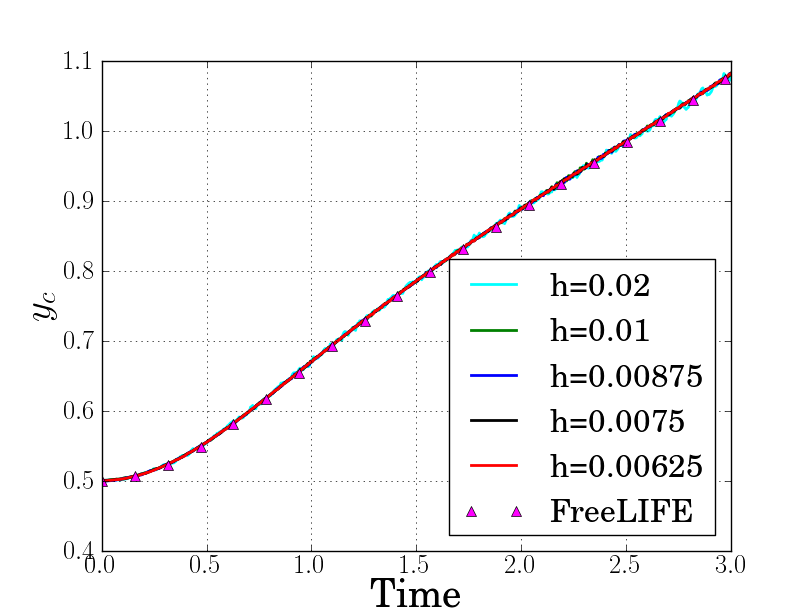}
    \label{subfig:elli_yc}
  }\\
  \subfigure[Vertical velocity.]{
    \includegraphics[width=0.4\textwidth]{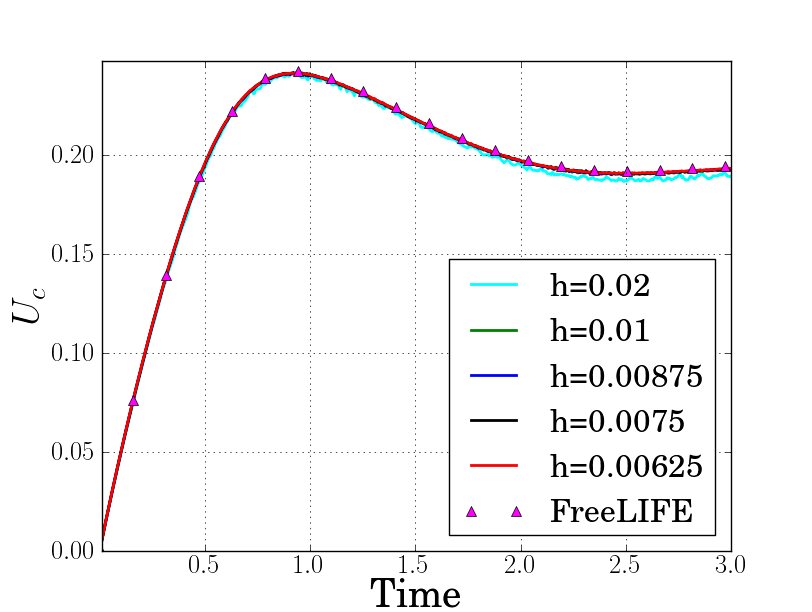}
    \label{subfig:elli_uc}
  }
  \subfigure[Circularity.]{
    \includegraphics[width=0.4\textwidth]{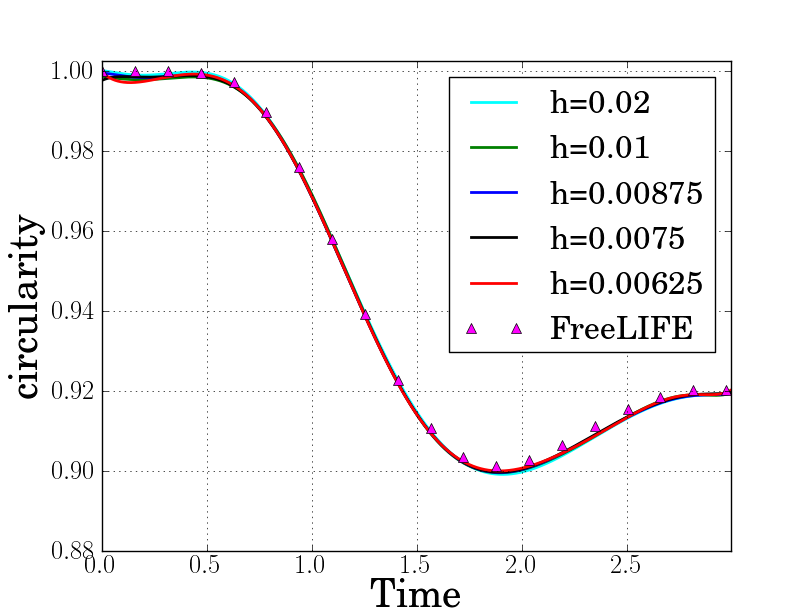}
    \label{subfig:elli_ci}
  }
  \caption{Results for the ellipsoidal bubble}
  \label{Fig:Ellipsoid}
\end{figure}

We run the simulations looking for solutions in finite element spaces spanned by
Lagrange polynomials of order $(2,1,1)$ for respectively the velocity, the
pressure and the level set.  In the first test case, the bubble reaches a
stationary circularity and its topology does not change. The velocity increases
until it attains a maximum then decreases to a constant
value. Figure~\ref{Fig:Ellipsoid} shows the results we obtained with different
mesh sizes. Three different groups presented their results in \cite{Hysing2009}
(\texttt{FreeLIFE}, \texttt{TP2D}, \texttt{MooNMD}). For the sake of clarity, we
only add on our graphs the data from one of the groups
(\texttt{FreeLIFE}). Nevertheless, the table~\ref{tab:elli_comp_res} shows a
comparison of our results with all groups published in \cite{Hysing2009}. We
monitor $c_{min}$ the minimum of the circularity, $t_{c_{min}}$ the time to
attain this minimum, $u_{c_{max}}$ the maximum velocity, $t_{u_{c_{max}}}$ the
time to reach it, and $y_c(t=3)$ the position of the bubble at final time
($t=3$).

\begin{figure}[h!tbp]%
  \centering
  \subfigure[Shape at final time ($t=3$).]{
    \includegraphics[width=0.35\textwidth]{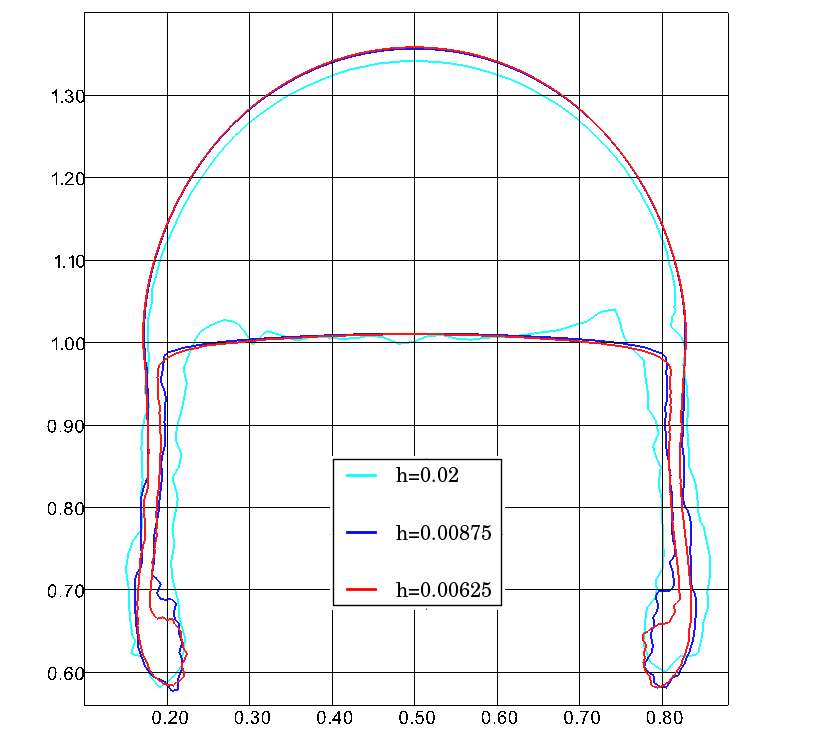}
    \label{subfig:ski_XY_paraview}
  }
  \subfigure[$y_c$ vertical.]{
    \includegraphics[width=0.4\textwidth]{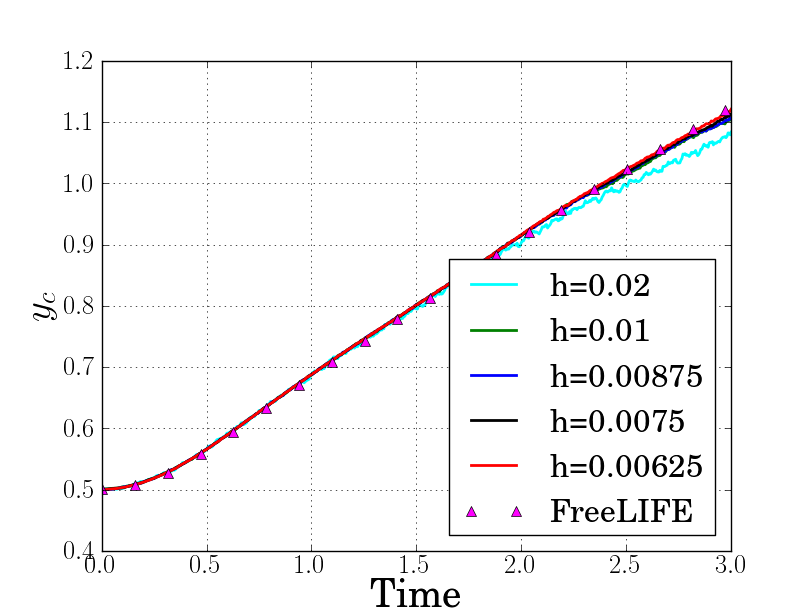}
    \label{subfig:skir_py}
  }\\
  \subfigure[Vertical velocity.]{
    \includegraphics[width=0.4\textwidth]{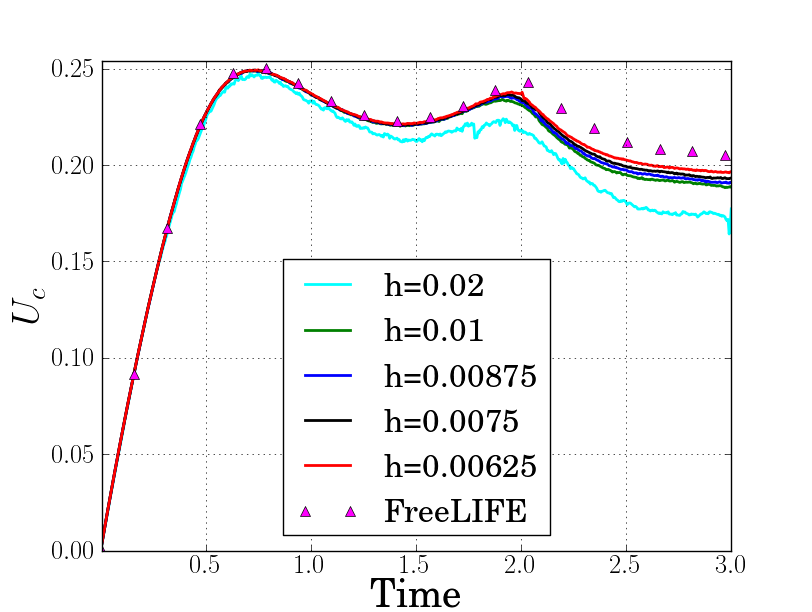}
    \label{subfig:skir_vy}
  }
  \subfigure[Circularity.]{
    \includegraphics[width=0.4\textwidth]{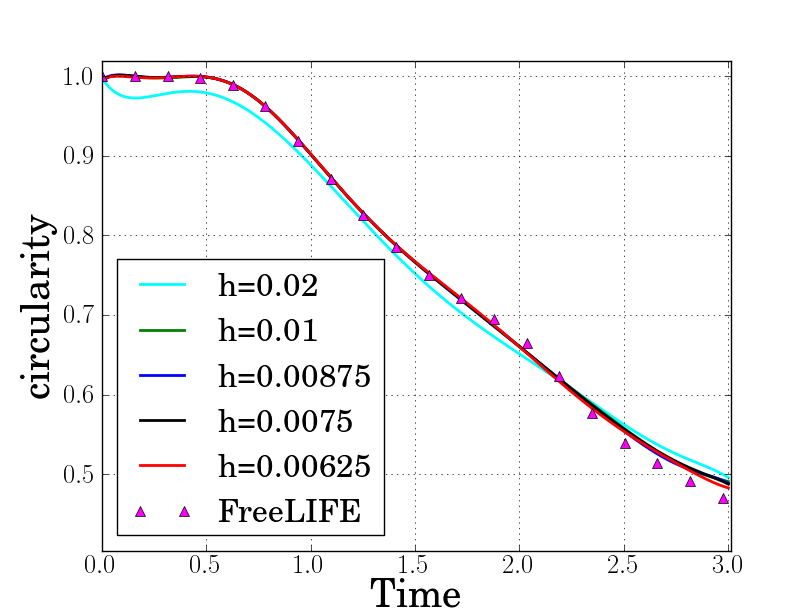}
    \label{subfig:skir_ci}
  }
  \caption{Results for the skirted bubble.}
  \label{Fig:Skirted}
\end{figure}

\begin{table}[htb!]
\normalsize
  \begin{center}
    \begin{tabular}{cccccc}
      \hline
               & $c_{min}$ & $t_{c_{min}}$  & $u_{c_{max}}$  & $t_{u_{c_{max}}}$ & $y_c(t=3)$   \\ \hline
   lower bound & 0.9011  & 1.8750  & 0.2417 & 0.9213  & 1.0799  \\
   upper bound & 0.9013  & 1.9041 & 0.2421  & 0.9313 &  1.0817  \\\hline \hline
   h=0.00625 & 0.9001  &1.9 &  0.2412&  0.9248 & 1.0815 \\
   h=0.0075 & 0.9001 & 1.9 & 0.2412 & 0.9251   & 1.0812 \\
   h=0.00875 &0.89998 &  1.9 &  0.2410&  0.9259& 1.0814  \\
   h=0.01 & 0.8999 & 1.9 & 0.2410 & 0.9252  & 1.0812 \\
   h=0.02 &0.8981  & 1.925 & 0.2400 & 0.9280 & 1.0787  \\
   \hline
    \end{tabular}
    \caption{Results comparison between benchmarks values (lower and upper bounds) and ours for ellipsoidal bubble.}
    \label{tab:elli_comp_res}
  \end{center}
\end{table}
In the second test case, the bubble gets more deformed because of the lower surface
tension. Some filaments (skirts) appear at the bottom. The velocity attains two
local maximum. Figure~\ref{Fig:Skirted} displays these
results and table~\ref{tab:ski_comp_res} shows the comparison with the benchmark
results. We monitor the same quantities as in the previous test case except that
we add the second maximum velocity $u_{c_{max_2}}$, and the time to reach it $t_{u_{c_{max_2}}}$.

\begin{table}[htb!]
  \normalsize
  \begin{center}
    \begin{tabular}{cccccccc}
      \hline
      & $c_{min}$ & $t_{c_{min}}$  & $u_{c_{max_1}}$  & $t_{u_{c_{max_1}}}$ &  $u_{c_{max_2}}$  & $t_{u_{c_{max_2}}}$ & $y_c(t=3)$   \\ \hline
      lower bound & 0.4647  & 2.4004  & 0.2502  & 0.7281 & 0.2393  & 1.9844 & 1.1249  \\
      upper bound & 0.5869   &3.0000  & 0.2524  & 0.7332 & 0.2440  & 2.0705 &  1.1380  \\\hline \hline
      h=0.00625&0.4616  &  2.995 & 0.2496 & 0.7574 & 0.2341 & 1.8828 & 1.1186 \\
      h=0.0075 & 0.4646 & 2.995 & 0.2495 & 0.7574 & 0.2333 & 1.8739 & 1.1111 \\
      h=0.00875& 0.4629 &  2.995&  0.2494 & 0.7565 & 0.2324&  1.8622&  1.1047 \\
      h=0.01 & 0.4642 & 2.995 & 0.2493 & 0.7559 & 0.2315 & 1.8522 & 1.1012 \\
      h=0.02 & 0.4744 & 2.995 & 0.2464 & 0.7529 & 0.2207 & 1.8319 & 1.0810 \\
      \hline
    \end{tabular}
    \caption{Results comparison between benchmarks values (lower and upper bounds) and ours for skirted bubble}
    \label{tab:ski_comp_res}
  \end{center}
\end{table}
Both tests show good agreements between our simulations and the ones from the
benchmark. We can notice that the final shape of the skirted bubble is very
sensitive to the mesh size and none of the groups agree on the exact shape which
can explain the differences that we see on the parameters in figure
\ref{Fig:Skirted} at time $t>2$.

\section{Vesicle dynamics simulation}
\label{sec:vesicle_dynamics_simulation}

\subsection{Model}
\label{sec:model}

The model usually admitted for vesicle membrane assumes three properties :
\textit{(i)} the membrane has a bending energy $E_b$, called the Canham,
Helfrich energy \cite{Canham1970, Helfrich1973}, \textit{(ii)} the inner fluid
is incompressible, so the total surface of the vesicle is conserved and finally
\textit{(iii)} the membrane is quasi inextensible, so the local perimeter is
conserved over the time.

\subsubsection{The bending energy}
\label{sec:bending-energy}
It has been shown \cite{Canham1970, Helfrich1973} that the bending energy is proportional to the square of the curvature of the membrane, in 2D it reads :
\begin{equation}
  E_b =  \int_{\Gamma} \frac{k_B}{2} \kappa^2,
  \label{eq:helfrich}
\end{equation}
where $k_B$ is the bending modulus (a typical value for phospholipidic membrane
is $k_B \approx 10^{-19}$ J).
Using the virtual power methods, the authors in \cite{Maitre2010} found a general expression (2D and 3D) for the force associated to this energy. It is given by :
\begin{equation}
  \bm{F_b} = \int_{\Omega} k_B \nabla \cdot \left[\frac{-\kappa^2}{2} \frac{\nabla \phi}{|\nabla \phi|} + \frac{1}{|\nabla \phi|} \left( \mathbb{I} - \frac{\nabla \phi \otimes \nabla \phi}{|\nabla \phi|^2} \right) \nabla \{|\nabla \phi| \kappa \} \right] \delta_{\varepsilon}.
\end{equation}

\subsubsection{Membrane inextensibility}
\label{sec:membr-inext}
The membrane inextensibility is equivalent numerically to impose that the
surface divergence of the velocity vanishes on the membrane. Several methods has
been developed to impose this constraint using Lagrange multipliers. This idea
has been applied with several different methods such as phase field method
\cite{Jamet2007, Maitre2010}, boundary integral methods
\cite{Veerapaneni20092334}, and even for level set methods \cite{Salac2011,
Laadhari2011}.  In particular, in \cite{Maitre2010}, the authors used a phase
field method in which the \textit{tension} is a variable defined in the entire
domain and it is given by the solution of an advection equation. This tension is
then added to the right hand side of the fluid equations as a force acting on
the membrane. The method from \cite{Salac2011}, is somehow similar to the
previous one in the sense that a \textit{tension} parameter is defined and also
added to a membrane force. Then, the system Navier-Stokes with tension equation
is solved by a 4 steps projection method discretized by finite differences.
In~\cite{Laadhari2011}, the authors used a level set method solved by
\textbf{FEM} in which the Lagrange multiplier is added to the variational
formulation of the Navier Stokes equations and acts as a pressure on the
membrane to keep the local perimeter constant. Moreover, two Lagrange
multipliers are added to the advection equation of the level set in order to
maintain constant the surface and the perimeter of the inner fluid. The mesh is
finally adaptively refined at each time step around the interface to get
improved accuracy.

In the present work, we use a method similar to the one presented
in~\cite{Laadhari2011}. Indeed, we add a Lagrange multiplier to the fluid
equations to impose the constraint \eqref{eq:stokes_3}. But our strategy is to
avoid adding Lagrange multipliers on advection equation (which are more difficult
to justify physically) and we choose to not adapt the mesh. Instead, we increase the
discretization orders to improve the perimeter conservation of the
vesicle.

Most of known behaviors of vesicles under flow take place at very low
Reynolds number. Thus, we assume that the fluid flow is governed by Stokes
equations subject to the inextensibility constraint of the membrane. We have :
  \begin{eqnarray}
    - 2 \nu D(\bm{u}) + \nabla p &=& \bm{F} \: \text{ in } \Omega \label{eq:stokes_1}\\
    \nabla \cdot \bm{u} &=& 0 \: \text{ in } \Omega \label{eq:stokes_2}\\
    \nabla_s \cdot \bm{u} &=& 0 \: \text{ on } \Gamma \label{eq:stokes_3}\\
    \bm{u} &=& \bm{g} \: \text{ on } \partial \Omega \label{eq:stokes_4}
  \end{eqnarray}
  with $\bm{u}$ the fluid velocity, $D(\bm{u}) = \dfrac{\nabla \bm{u} + \nabla
    \bm{u}^T}{2}$ the deformation tensor, $p$ the pressure, $\bm{F}$ the
  external forces, and $\nabla_s \cdot \bm{u} = \nabla \cdot \bm{u} - (\nabla
  \bm{u} \cdot \bm{n} ) \cdot \bm{n}$ the surfacic divergence. The variational
  formulation associated to the problem
  \eqref{eq:stokes_1}-\eqref{eq:stokes_2}-\eqref{eq:stokes_3}-\eqref{eq:stokes_4}
  reads: Find $(\bm{u}, p, \lambda) \in V \times L^2_0(\Omega) \times
  H^{1/2}(\Gamma)$ which verify $\forall (\bm{v}, q, \mu) \in H^1_0(\Omega)^2
  \times L^2_0(\Omega) \times H^{1/2}(\Gamma)$ :
  \begin{eqnarray}
    2 \int_{\Omega} \nu(\phi) D(\bm{u}) : D(\bm{v}) - \int_{\Omega} p \nabla \cdot \bm{v} + \int_{\Gamma} \lambda \nabla_s \cdot \bm{v} &=& \int_{\Omega} \bm{F}(\phi) \cdot \bm{v}, \label{eq:varstokes_1}\\
    \int_{\Omega} q \nabla \cdot \bm{u} &=& 0, \label{eq:varstokes_2} \\
    \int_{\Gamma} \mu \nabla_s \cdot \bm{u} &=& 0, \label{eq:varstokes_3}
  \end{eqnarray}
  with $\lambda$ the Lagrange multiplier associated to the free surfacic
  divergence and $V = \{ \bm{v} \in H^1(\Omega)^2 \ | \ \bm{v}|_{\partial
  \Omega} = \bm{g} \} $.
  It is obvious in this formulation that the Lagrange multiplier can be
  interpreted as a pressure acting on the interface.

  The problem
  \eqref{eq:varstokes_1}-\eqref{eq:varstokes_2}-\eqref{eq:varstokes_3} is
  coupled with the level set advection \eqref{eq:varadvection} and discretized
  using the finite element method. Thereby we introduce ${\mathbb U}_h^{n+1}$,
  ${\mathbb P}_h^n$, ${\mathbb Q}_h^m$ ans ${\mathbb R}_h^k$ the discrete finite
  elements spaces depending on mesh size $h$ and based on Lagrange polynomials
  of degree $n+1$, $n$, $m$ and $k$ for the velocity, pressure, Lagrange
  multipliers and levelset respectively. A complete description of the strategy
  to obtain high order level set method (including reinitialization at high
  order and benchmarks) will be presented in \cite{doyeux_ls}.  Let $(\bm{u}_h,
  p_h, \lambda_h, \phi_h) \in {\mathbb U}_h^{n+1} \times {\mathbb P}_h^n \times
  {\mathbb Q}_h^m \times {\mathbb R}_h^k$ be the discretization of $(\bm{u}, p,
  \lambda, \phi)$. The discrete version of
  \eqref{eq:varstokes_1}-\eqref{eq:varstokes_2}-\eqref{eq:varstokes_3}-\eqref{eq:varadvection}
  reads: Find $(\bm{u}_h, p_h, \lambda_h, \phi_h) \in {\mathbb U}_h^{n+1} \times
  {\mathbb P}_h^n \times {\mathbb Q}_h^m \times {\mathbb R}_h^k$ which verify
  $\forall (\bm{v}_h, q_h, \mu_h, \psi_h) \in {\mathbb U}_h^{n+1} \times
  {\mathbb P}_h^n \times {\mathbb Q}_h^m \times {\mathbb R}_h^k$ :
  \begin{eqnarray}
    2 \int_{\Omega} \nu(\phi_h) D(\bm{u}_h) : D(\bm{v}_h) - \int_{\Omega} p_h \nabla \cdot \bm{v}_h + \int_{\Omega} \lambda_h \nabla_s \cdot \bm{v}_h \  \delta_{\varepsilon}(\frac{\phi_h}{|\nabla \phi_h|}) &=& \int_{\Omega} \bm{F}_h \cdot \bm{v}_h, \label{eq:discr_varstokes_1}\\
    \int_{\Omega} q_h \nabla \cdot \bm{u}_h &=& 0, \label{eq:discr_varstokes_2} \\
    \int_{\Omega} \mu_h \nabla_s \cdot \bm{u}_h \ \delta_{\varepsilon}(\frac{\phi_h}{|\nabla \phi_h|}) &=& 0, \label{eq:discr_varstokes_3}\\
    \int_{\Omega} \frac{\partial \phi_h}{\partial t}  \psi_h + \int_{\Omega} (\bm{u}_h \cdot \nabla \phi_h) \psi_h + \int_{\Omega} S(\phi_h, \psi_h) &=& 0. \label{eq:discr_varadvection}
  \end{eqnarray}

  Note that we replaced the integral over $\Gamma$ by integral over $\Omega$
  thanks to the delta function. Thus, the Lagrange multiplier space is defined
  only in the elements in a region of $2 \varepsilon$ around the interface $\phi
  = 0$. In practice, from the implementation point of view, we only add to the
  global matrix the coefficients that correspond to these ``few'' elements.

\subsection{Tank treading motion}

The first behavior on which we validate our model is the tank treading motion
(\textbf{TT}). In a linear shear flow, if the viscosity ratio between the inner and
outer fluids $ \nu_r = \dfrac{\nu_2}{\nu_1}$ is lower than a critical value, the
vesicle reaches a steady angle with respect to the horizontal. At the same time,
the membrane is rotating along the vesicle with a constant velocity (this motion
is similar to the chain of a tank hence the name \textit{tank treading motion}).

We need to define dimensionless parameters which control this system: \textit{(i)} the
reduced area $\alpha$ is the ratio between the area of the vesicle ($A$) and the
area of a circle having the same perimeter ($P$) $\bigg(\alpha = \dfrac{4 \pi
A}{P^2}\bigg)$, \textit{(ii)} the capillary number $C_a$ which is the ratio between the characteristic
time of the shear $(\frac{1}{\gamma})$ and a characteristic time related to the
curvature force, it is given by $C_a = \dfrac{\nu_2 \gamma R_0^3}{k_B}$ with
$R_0=\dfrac{P}{2\pi}$ stands for a typical size of the vesicle. For high values of $C_a$,
the particle is more deformable  due to the hydrodynamic forces, while for low
$C_a$, the Helfrich energy requires a higher cost  to change the curvature of
the vesicle.  The confinement $c$ is defined as the ratio between the equivalent
radius of the vesicle and the half width of the channel : $c=\dfrac{R_0}{L/2}$.

\begin{figure}[h!tbp]%
  \center
  \subfigure[$t = 0$]{
    \includegraphics[width=0.48\linewidth]{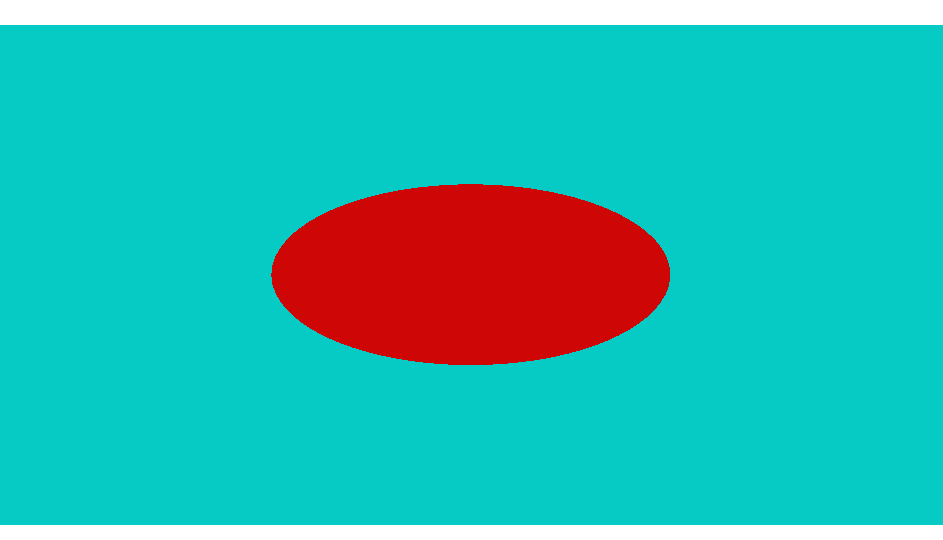}
    \label{subfig:vesicle_tanktreading_init}
  }
  \subfigure[$t = 5$]{
    \includegraphics[width=0.48\linewidth]{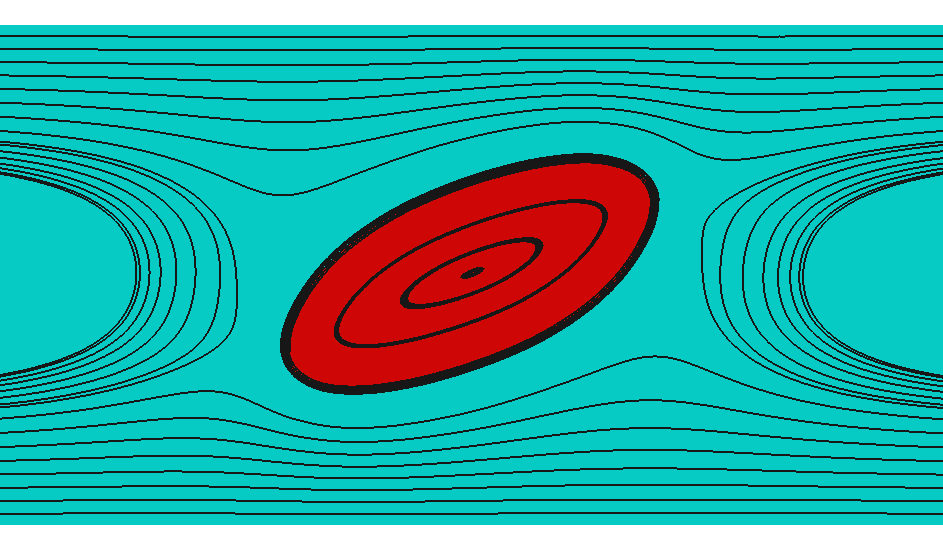}
    \label{subfig:vesicle_tanktreading_steady}
  }
  \caption{Vesicle at $\alpha=0.8, C_a = 0.12 \text{ and } c = 0.6$ reaches a steady angle. Streamlines are along the membrane, the vesicle is under tank treading regime.}
  \label{Fig:vesicle_tanktreading}
\end{figure}

Figure \ref{Fig:vesicle_tanktreading} shows the initial and steady states of a
vesicle in a shear flow with $\alpha = 0.8$, $C_a=0.12$ and $c=0.6$. Initially,
the vesicle is placed horizontally and initialized as an ellipse. At steady state, the vesicle
has taken a steady angle with respect to the horizontal, and the fluid velocity
rotates along the membrane (this is shown by streamlines in figure
\ref{subfig:vesicle_tanktreading_steady}). The simulation has been run in a
rectangular box of size $[10, 2.12]^2$ descretized by $12400$ elements. The time
step was taken as $\delta t = 5 \times 10^{-3}$ and the finite elements
discretization space was ${\mathbb U}_h^{3} \times {\mathbb P}_h^2 \times
{\mathbb Q}_h^2 \times {\mathbb R}_h^2$.

\begin{figure}[h!tbp]%
  \center
  \subfigure[Vesicle angle]{
    \includegraphics[width=0.48\linewidth]{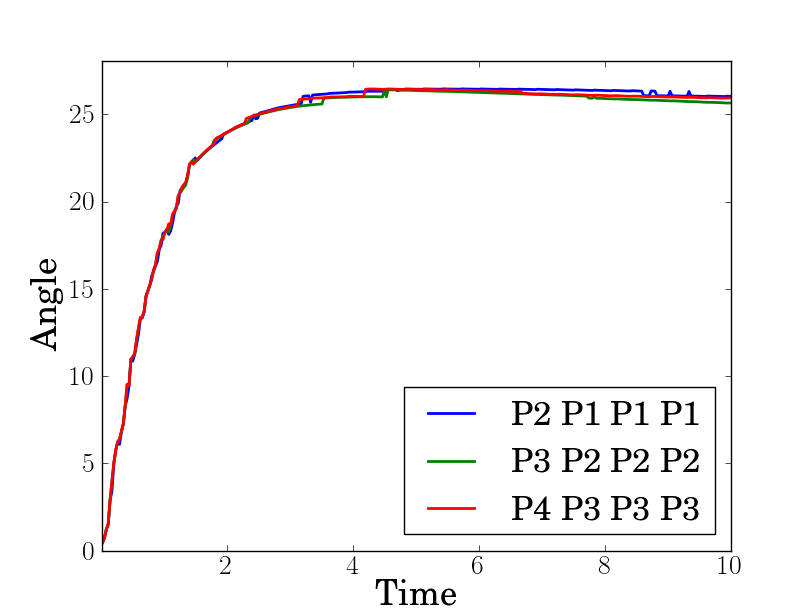}
    \label{subfig:angle_tt_order}
  }
  \subfigure[Loss of perimeter (in \%)]{
    \includegraphics[width=0.48\linewidth]{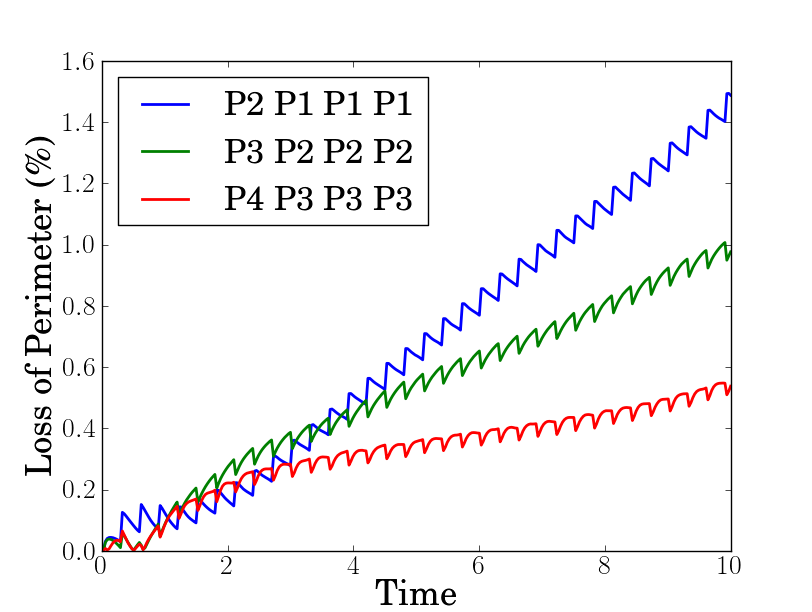}
    \label{subfig:perim_tt_order}
  }
  \caption{Vesicle angle and loss of perimeter for different polynomial approximation order. The legend gives in this order : $(n+1,n,m,k)$.}
  \label{Fig:tt_order}
\end{figure}

Different finite element discretization has been tested. Figure
\ref{Fig:tt_order} shows the results of three different sets of polynomials
orders $(n,m,k)$ : $(1,1,1)$, $(2,2,2)$ and $(3,3,3)$. The computational time is
multiplied by $2$ from $(n,m,k) = (1,1,1)$ to $(n,m,k) = (2,2,2)$ and multiplied
by $4$ from $(n,m,k) = (1,1,1)$ to $(n,m,k) = (3,3,3)$. We can see on figure
\ref{subfig:angle_tt_order} that the tank treading steady angle doesn't change
dramatically for simulations up to this final time. Nevertheless, the loss of
perimeter is a crucial point for accuracy in long time simulations. Thus, we
plotted in figure \ref{subfig:perim_tt_order} the loss of perimeter :
$\dfrac{|p-p_0|}{p_0} \times 100$ for the different polynomial approximations
used. The oscillations seen on figure \ref{subfig:perim_tt_order} are due to the
reinitialization steps. One can see that increasing the polynomial order
approximation improves the conservation of the perimeter. The exact role of
each polynomial approximation on the perimeter conservation accuracy still has
to be investigated. For the other simulations, we choose to take $(n,m,k) =
(2,2,2)$ as polynomial order approximation which, for our applications, seems to
be a good compromise between accuracy and computational time.

\begin{figure}[h!tbp]%
  \centering
  \subfigure[Vesicle angle in shear flow for $c=0.6$ and $C_a = 0.1$]{
  \includegraphics[width=0.45\linewidth]{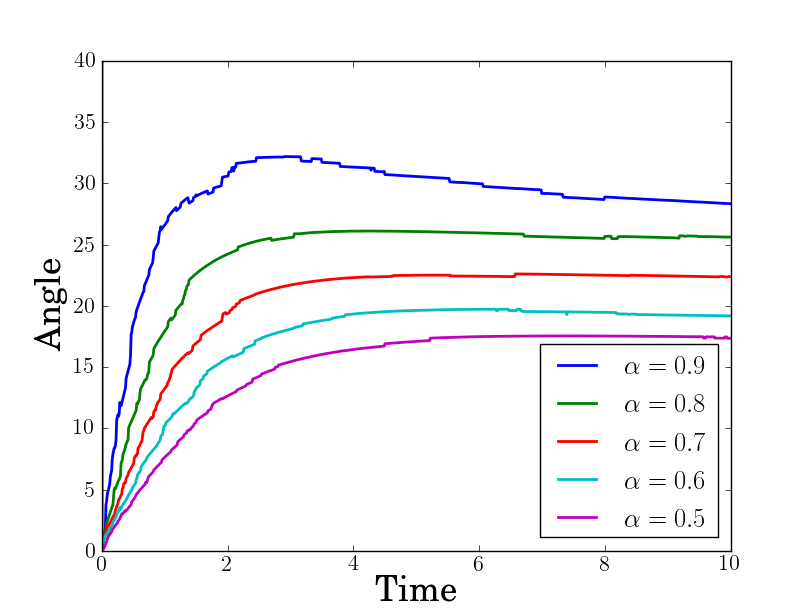}
  \label{subfig:vesicle_angles_rv}
  }
  \subfigure[Steady angle of vesicle in tank treading motion as a function of reduced area for different confinements.]{
  \includegraphics[width=0.45\linewidth]{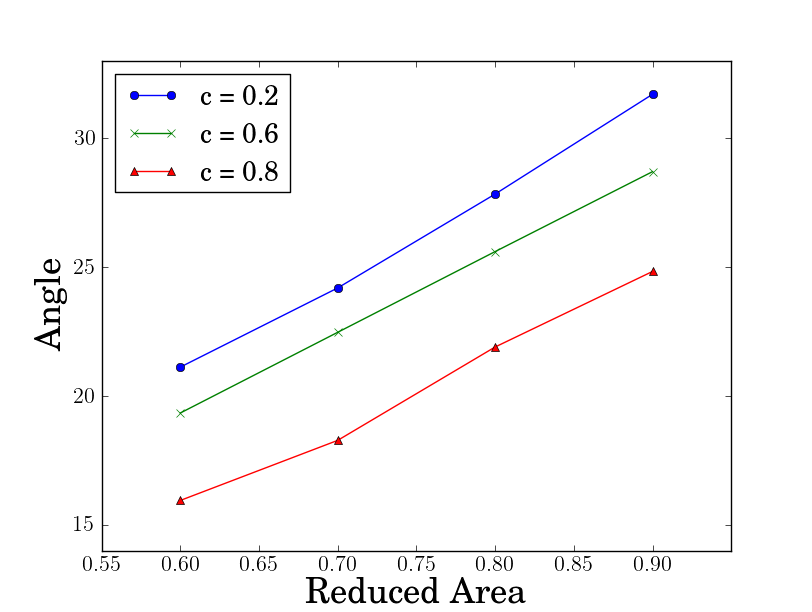}
  \label{subfig:vesicle_angles_conf}
  }
  \caption{Vesicle angle as a function of time, reduced area and confinement.}
  \label{Fig:vesicle_angles}
\end{figure}

It has been shown in \cite{Kaoui2011} that the steady tank treading angle decreases
with the reduced area. So, we run the same simulation than shown in
\ref{Fig:vesicle_tanktreading} by changing the reduced area of the vesicle. The
steady angles exhibit the expected behavior as one can see it in figure
\ref{subfig:vesicle_angles_rv}.

Moreover, it has been also shown in \cite{Kaoui2011} that for a given reduced area, the
steady angle is lower for high confinements. Figure
\ref{subfig:vesicle_angles_conf} shows the steady angles that we found for three
different confinements at different reduced volumes. Once again, we obtain the
expected behavior.

\subsection{Tumbling motion}

When the viscosity ratio $\nu_r$ is above a critical value, the vesicle tends to
follow a solid rotation, it is called the \textit{tumbling motion}. Increasing the
viscosity ratio increases the rotation frequency of the vesicle.

\begin{figure}[h!tbp]%
  \center
  \subfigure[$t = 1.0$]{
    \includegraphics[width=0.3\linewidth]{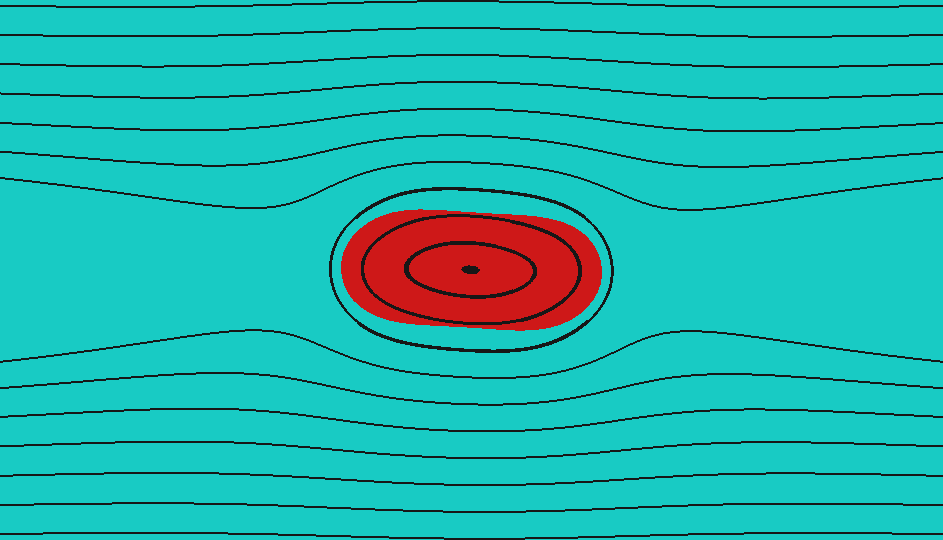}
    \label{subfig:}
  }
  \subfigure[$t = 1.1$]{
    \includegraphics[width=0.3\linewidth]{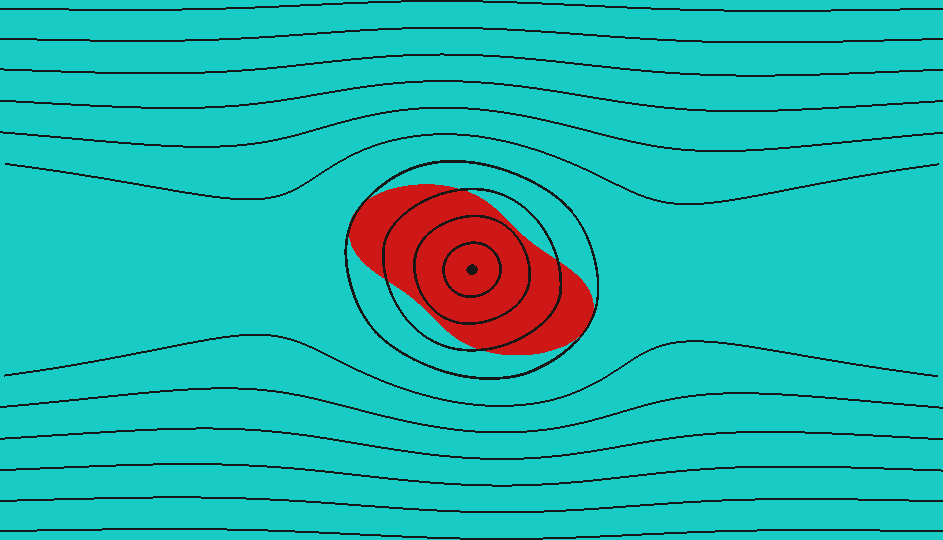}
    \label{subfig:}
  }
  \subfigure[$t=6.0$]{
    \includegraphics[width=0.3\linewidth]{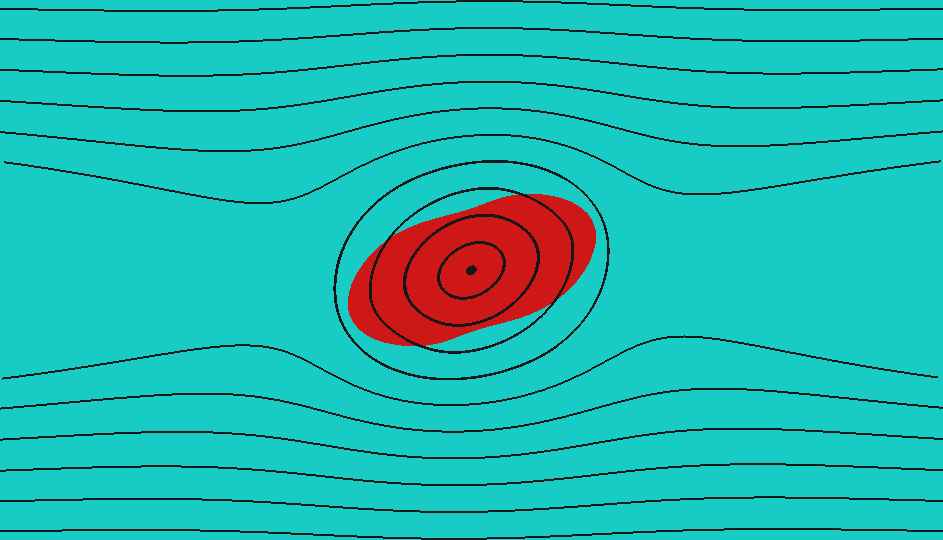}
    \label{subfig:}
  }
  \caption{Tumbling of a vesicle}
  \label{Fig:Tumbling}
\end{figure}

To reproduce this behavior, we set initially a vesicle as an ellipse in a box of
size $[14,5]^2$, discretized with $7218$ elements. We chose a time step of
$\delta t = 3 \times 10^{-2}$. These parameters are not as refined than in the
previous section in order to get long time simulations in a reasonable
computational time. The other parameters of the simulation were : $C_a = 6
\times 10^{-2}$, $c=0.25$ and $\alpha = 0.8$. Figure \ref{Fig:Tumbling} displays
the tumbling motion.

\begin{figure}[h!tbp]%
  \centering
  \subfigure[Tumbling angle as a function of time for two different viscosity ratios.]{
  \includegraphics[width=0.48\linewidth]{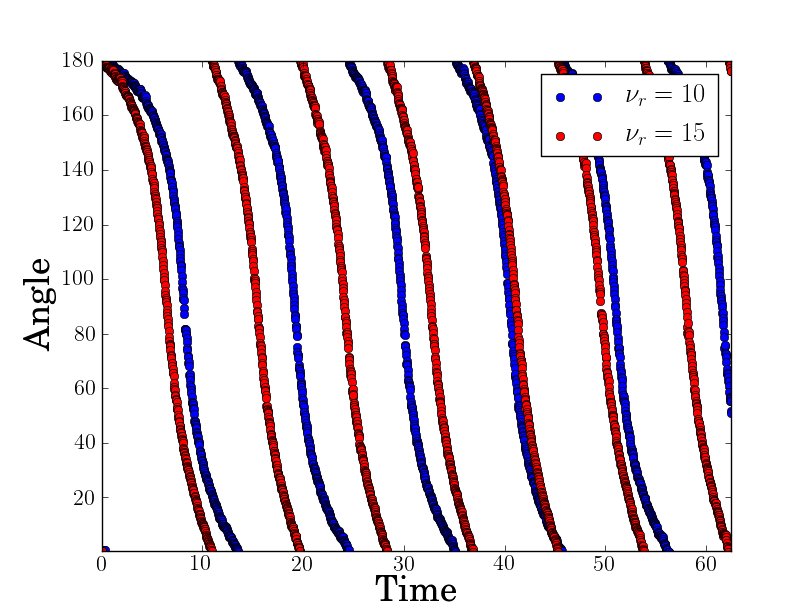}
  \label{subfig:tumbling_angle}
  }
  \subfigure[Tumbling frequency for different viscosity ratios.]{
  \includegraphics[width=0.48\linewidth]{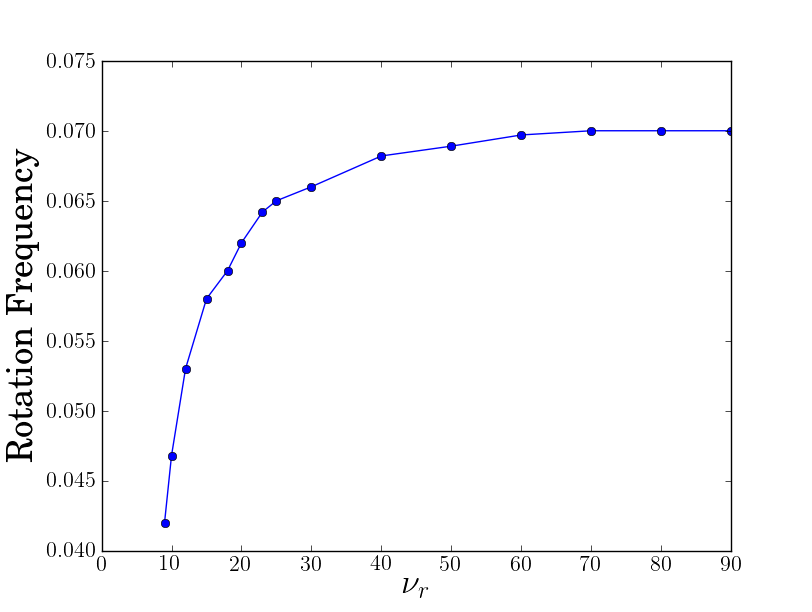}
  \label{subfig:tumbling_frequency}
  }
  \caption{Tumbling angle and frequency.}
  \label{Fig:Tumbling_angles_frequency}
\end{figure}

Increasing the viscosity ratio increases the rotation frequency of the tumbling
as it can be seen on figures~\ref{subfig:tumbling_angle}
and~\ref{subfig:tumbling_frequency}. Moreover, as described in
\cite{ghigliotti2010}, when one increases the viscosity ratio, the rotation
frequency of the vesicle reaches a steady value which corresponds to the steady rotation of
a solid object in a shear flow. We can see this phenomenon in
figure~\ref{subfig:tumbling_frequency}.

\section*{Conclusion}

We have presented a new numerical framework for the simulation of vesicle under
flow. This framework is based on level set methods solved by a (possibly high
order) finite element method. First the level set framework for two-fluid flows
at order $(1,1,1)$ has been verified using a numerical benchmark. Then, for more
complex entities such us vesicles, the incompressibility of the membrane is
taken into account using a Lagrange multiplier defined on the interface.
Compared to the literature, one of the novelties of our work lies in the
possibility of taking into account this constraint without refining the mesh or
introducing auxiliary variables. This being done at a lower cost in the sense
that this multiplier affects only the elements crossed by the interface.

A validation on tank treading and tumbling regimes has been presented. In
particular, for the tank treading motion we have used high order polynomial
approximations $(2, 2, 2)$ and $(3, 3, 3)$. 
However, a theoretical study of the role of each
approximation order and its verification on numerical tests has to be done in
a near future. Finally other basic behaviors of vesicles under flow has been tested
for a polynomial order of $(2,2,2)$ and show the expected results in good
agreement with the literature.

We are also interested in clustering phenomenon which is very important in the case of red
blood cells and also to do 3D calculations (the \texttt{FEEL++} library allows
it without much changes in the code).

\section*{Acknowledgments}

The authors would like to thank the R\'egion Rh\^one-Alpes (ISLE/CHPID project)
as well as the French National Research Agency (the MOSICOB and Cosinus-HAMM
projects) for their financial support.

\section*{References}

\bibliographystyle{elsarticle-num}
\bibliography{ls}

\begin{thebibliography}{10}
\expandafter\ifx\csname url\endcsname\relax
  \def\url#1{\texttt{#1}}\fi
\expandafter\ifx\csname urlprefix\endcsname\relax\def\urlprefix{URL }\fi
\expandafter\ifx\csname href\endcsname\relax
  \def\href#1#2{#2} \def\path#1{#1}\fi

\bibitem{Kaoui2011}
B.~Kaoui, J.~{Harting}, C.~{Misbah}, {Two-dimensional vesicle dynamics under
  shear flow: Effect of confinement}, pre 83~(6) (2011) 066319--+.
\newblock \href {http://arxiv.org/abs/1011.6061} {\path{arXiv:1011.6061}},
  \href {http://dx.doi.org/10.1103/PhysRevE.83.066319}
  {\path{doi:10.1103/PhysRevE.83.066319}}.

\bibitem{Beaucourt2004}
J.~Beaucourt, F.~Rioual, T.~Seacuteon, T.~Biben, C.~Misbah, Steady to unsteady
  dynamics of a vesicle in a flow, Phys. Rev. E 69~(1) (2004) 011906--.

\bibitem{Maitre2010}
E.~Maitre, C.~Misbah, P.~Peyla, A.~Raoult, Comparison between advected-field
  and level-set methods in the study of vesicle dynamics, ArXiv e-prints\href
  {http://arxiv.org/abs/1005.4120} {\path{arXiv:1005.4120}}.

\bibitem{Salac2011}
D.~Salac, M.~Miksis, A level set projection model of lipid vesicles in general
  flows, Journal of Computational Physics 230~(22) (2011) 8192--8215.

\bibitem{Laadhari2011}
A.~Laadhari, P.~Saramito, C.~Misbah,
  \href{http://hal.archives-ouvertes.fr/hal-00604145/en/}{{Computing the
  dynamics of biomembranes by combining conservative level set and adaptive
  finite element methods}}, cNRS (Jun. 2011).
\newline\urlprefix\url{http://hal.archives-ouvertes.fr/hal-00604145/en/}

\bibitem{ismail_lefebvre:2012}
M.~Ismail, A.~Lefebvre, A ``necklace'' model for vesicles simulation in 2d,
  submitted to Journal of Computational Physics (2012).

\bibitem{Winkelmann2007}
C.~Winkelmann, Interior penalty finite element approximation of navier-stokes
  equations and application to free surface flows, Ph.D. thesis (2007).

\bibitem{Osher1988}
S.~Osher, J.~A. Sethian, Fronts propagating with curvature dependent speed:
  Algorithms based on hamilton-jacobi formulations, Journal of computational
  physics 79~(1) (1988) 12--49.

\bibitem{book_Sethian}
J.~Sethian, Level Set Methods and Fast Marching Methods, Cambridge University
  Press, 1996.

\bibitem{book_Osher}
R.~F. Stanley~Osher, Level Set Methods and Dynamic Implicit Surfaces, Springer,
  S.S. Antman, J.E. Marsden, L. Sirovich.

\bibitem{Cottet2010}
E.~M. Georges-Henri~Cottet, A level set method for fluid-structure interactions
  with immersed surfaces, Mathematical Models and Methods in Applied Sciences.

\bibitem{PRUDHOMME:2012:HAL-00662868:3}
C.~Prud'Homme, V.~Chabannes, V.~Doyeux, M.~Ismail, A.~Samake, G.~Pena,
  \href{http://hal.archives-ouvertes.fr/hal-00662868}{{Feel++: A Computational
  Framework for Galerkin Methods and Advanced Numerical Methods}}, submitted to
  ESAIM Proc. (Jan. 2012).
\newline\urlprefix\url{http://hal.archives-ouvertes.fr/hal-00662868}

\bibitem{prudhomme:_feel}
C.~Prud'homme, V.~Chabannes, G.~Pena, Feel++: {F}inite {E}lement {E}mbedded
  {L}anguage in {C}++, Free Software available at http://www.feelpp.org,
  contributions from A. Samake, V. Doyeux, M. Ismail and S. Veys.

\bibitem{Prudhomme.ea:2006}
C.~Prud'homme, A domain specific embedded language in {C++} for automatic
  differentiation, projection, integration and variational formulations,
  Scientific Programming 14.

\bibitem{Franca1992}
L.~P. Franca, S.~L. Frey, T.~J. Hughes, Stabilized finite element methods: I.
  application to the advective-diffusive model, Computer Methods in Applied
  Mechanics and Engineering 95~(2) (1992) 253--276.

\bibitem{Burman2006}
E.~Burman, P.~Hansbo, Edge stabilization for the generalized stokes problem: A
  continuous interior penalty method, Computer Methods in Applied Mechanics and
  Engineering 195~(19-22) (2006) 2393--2410.

\bibitem{chabannes11:_high}
V.~Chabannes, G.~Pena, C.Prud'homme, High order fluid structure interaction in
  2d and 3d. application to blood flow in arteries, in: Fifth International
  Conference on Advanced COmputational Methods in ENgineering (ACOMEN 2011),
  2011.

\bibitem{Hysing2009}
S.~Hysing, S.~Turek, D.~Kuzmin, N.~Parolini, E.~Burman, S.~Ganesan, L.~Tobiska,
  Quantitative benchmark computations of two-dimensional bubble dynamics,
  International Journal for Numerical Methods in Fluids 60~(11) (2009)
  1259--1288.
\newblock \href {http://dx.doi.org/10.1002/fld.1934}
  {\path{doi:10.1002/fld.1934}}.

\bibitem{Canham1970}
P.~Canham, The minimum energy of bending as a possible explanation of the
  biconcave shape of the human red blood cell, Journal of Theoretical Biology
  26~(1) (1970) 61 -- 81.
\newblock \href {http://dx.doi.org/DOI: 10.1016/S0022-5193(70)80032-7}
  {\path{doi:DOI: 10.1016/S0022-5193(70)80032-7}}.

\bibitem{Helfrich1973}
W.~Helfrich, Elastic properties of lipid bilayers: theory and possible
  experiments., Z Naturforsch C 28~(11) (1973) 693--703--.

\bibitem{Jamet2007}
D.~Jamet, C.~Misbah, Towards a thermodynamically consistent picture of the
  phase-field model of vesicles: Local membrane incompressibility, Phys. Rev. E
  76~(5) (2007) 051907.
\newblock \href {http://dx.doi.org/10.1103/PhysRevE.76.051907}
  {\path{doi:10.1103/PhysRevE.76.051907}}.

\bibitem{Veerapaneni20092334}
S.~K. Veerapaneni, D.~Gueyffier, D.~Zorin, G.~Biros, A boundary integral method
  for simulating the dynamics of inextensible vesicles suspended in a viscous
  fluid in 2d, Journal of Computational Physics 228~(7) (2009) 2334 -- 2353.
\newblock \href {http://dx.doi.org/10.1016/j.jcp.2008.11.036}
  {\path{doi:10.1016/j.jcp.2008.11.036}}.

\bibitem{doyeux_ls}
V.~Doyeux, C.~Prud'homme, M.~Ismail, A framework toward high order level set
  method., in preparation.

\bibitem{ghigliotti2010}
G.~Ghigliotti, T.~Biben, C.~Misbah, Rheology of a dilute two-dimensional
  suspension of vesicles, Journal of Fluid Mechanics 653 (2010) 489--518.
\newblock \href {http://dx.doi.org/10.1017/S0022112010000431}
  {\path{doi:10.1017/S0022112010000431}}.

\end{thebibliography}

\end{document}